\documentclass[11pt,a4paper]{article}
\usepackage[utf8]{inputenc}
\usepackage[english]{babel}
\usepackage{amsmath}
\usepackage{amsfonts}
\usepackage{wrapfig}
\usepackage{graphicx}
\usepackage{booktabs}
\usepackage{amssymb}
\usepackage[left=3cm,right=3cm,top=2cm,bottom=2cm]{geometry}

\usepackage[utf8]{inputenc}
\usepackage[english]{babel}
\usepackage{microtype}
\usepackage{mathrsfs}
\usepackage{amsthm, thmtools}
\usepackage{amsfonts}
\usepackage{amsmath}
\usepackage{amssymb}
\usepackage{amscd}
\usepackage{mathtools}
\usepackage{wasysym}
\usepackage{braket}
\usepackage{hyperref}
\usepackage{appendix}
\usepackage{enumitem}
\usepackage{framed}
\usepackage[all]{xy}
\usepackage{xfrac}
\usepackage{color}
\usepackage{bbm}
\usepackage{stmaryrd}
\usepackage{fancyhdr}
\usepackage{calligra}
\usepackage{bigints}
\usepackage{afterpage}

\newtheorem{teorema}{Theorem}

\newtheorem{co}{Corollary}
\newtheorem{lemma}{Lemma}

\theoremstyle{definition}

\renewcommand{\epsilon}{\varepsilon}
\addto{\captionsenglish}{\renewcommand{\bibname}{ References}}

\author{Chiara Bellotti-Giuseppe Puglisi}
\title{On the Deuring-Heilbronn Phenomenon}
\begin{document}
\begin{center} 
		{\LARGE{\bf Elementary methods in the study of}}
		\vspace{2mm}
		
		{\LARGE{\bf Deuring-Heilbronn Phenomenon}}\\ 
		\vspace{2mm} 
		{\LARGE{\bf $ $ }}\\ 
		\vspace{5mm} 
		{\large Chiara Bellotti\footnote{University of Pisa, Department od Mathematics, chbellotti@gmail.com}, Giuseppe Puglisi\footnote{University of Pisa, Department of Mathematics, giuseppe.puglisi@unipi.it}}\\
		\vspace{2mm}
		\textit{University of Pisa, Pisa, Italy}
	\end{center} 
	\vspace{5mm} 
	\textbf{Abstract.} \textit{The aim of this work is to improve some elementary results regarding both the Deuring-Phenomenon and the Heilbronn-Phenomenon. We will give better estimates regarding both the influence of zeros of the Riemann zeta function  on the exceptional zeros and that of the non-trivial zeros of arbitrary $L$-functions belonging to non-principal characters on the exceptional zeros.} 
\section{Introduction}
Let $L(s,\chi_{D})$ be a Dirichlet $L$-function belonging to the real primitive character $\chi_{D}$ modulus $D$ satisfying $\chi_{D}(-1)=-1$. Let $h(-D)$ be the class number of the imaginary quadratic field $\mathbb{Q}(\sqrt{-D})$.

Two conjectures involving the class number $h(-D)$ of the imaginary quadratic field belonging to the fundamental discriminant $-D<0$ were raised by Gauss, who published them in 1801 \cite{gaussconj}. The first problem was about determining all the negative fundamental discriminants with class number one. The second problem was about proving that $h(-D) \rightarrow \infty$, as $D \rightarrow \infty$.

Regarding the second conjecture, in 1913, Gronwall \cite{gron} proved that if the function $L\left(s, \chi_{D}\right)$ belonging to the real primitive character $\chi_{D}(n)=\left(\frac{-D}{n}\right)$ has no zero in the interval
$\left[1-\frac{\alpha}{\log D}, 1\right]$, then
$h(-D)>\frac{b(\alpha) \sqrt{D}}{\log D \sqrt{\log \log D}}
$, where $\alpha$ is a constant and $b(\alpha)$ is a constant depending only on $\alpha$.

In 1918, Hecke \cite{landau} proved that, under the same hypotheses of Gronwall's theorem, the inequality $h(-D)>\frac{b^{\prime}(\alpha) \sqrt{D}}{\log D}$ holds, where $\alpha$ is a constant and $b^{\prime}(\alpha)$ is a constant depending only on $\alpha$.

In $1933$, Deuring \cite{dg} proved that under the assumption of the falsity of the classical Riemann Hypothesis the relation  $h(-D)\geq 2$ holds for $D>D_{0}$, where $D_{0}$ is a constant.\\In $1934$, Mordell \cite{mg} improved the result found by Deuring. Under the assumption of the falsity of the classical Riemann Hypothesis, Mordell proved that $h(-D)\rightarrow\infty$ as $D\rightarrow\infty$.\\These results showed an interesting connection between the possibly existing real zeros of special $L$-functions and the non-trivial zeros of the $\zeta$-function.

Better results regarding the influence of zeros of $\zeta(s)$ on the exceptional zeros, or equivalently, the \textit{Deuring phenomenon}, were provided by the work of Pintz, who used a new approach involving some elementary methods.\\ In $1976$, Pintz \cite{pintz3} proved that, assuming a relatively strong upper bound for $h(-D)$, it is possible to determine, up to a factor $1+o(1)$, the values of the corresponding $L$-function in a large domain of the critical strip. \medskip\\
\textbf{Theorem. (Pintz)} \textit{Given $0<\epsilon<1 / 8$ and $D>D_{1}(\epsilon)$, where $D_{1}(\epsilon)$ is an effective constant depending on $\epsilon$, we define the domain $H(\epsilon, D)$, depending on $\epsilon$ and on $D$, as the set
$$
\begin{array}{r}
H(\epsilon, D)=\left\{s ; s=1-\tau+i t,|1-s| \geq 1 / \log ^{4} D, 0 \leq \tau \leq \frac{1}{4}-\epsilon\right. \\
\left.|s| \leq D^{\left(\frac{1}{4}-\frac{\epsilon}{2}\right) \frac{1}{\varrho}-\frac{3}{4}} \text { where } \varrho=\max \left(\tau, D^{-\epsilon / 4}\right)\right\}
\end{array}
$$
If the inequality
$$
h(-D) \leq(\log D)^{3 / 4}
$$holds, then neither $L(s,\chi_{D})$ nor $\zeta(s)$ has a zero in $H(\epsilon, D)$, and for $s \in H(\epsilon, D)$, we have
$$
L(s,\chi_{D})=\frac{\zeta(2 s)}{\zeta(s)} \prod_{p \mid D}\left(1+\frac{1}{p^{s}}\right)\left[1+O\left(\exp \left\{-\frac{1}{8} (\log  D)^{1 / 4}\right\}\right)\right]
$$}
 An immediate consequence is that, except for the possible Siegel zero, neither $L(s, \chi_{D})$ nor $\zeta(s)$ has a zero in this domain. Also, a weakened form of Mordell's theorem follows, namely that if $h(-D) \not\rightarrow \infty$ for $D \rightarrow \infty$, then $\zeta(s)$ has no zero in the half-plane $\sigma>\frac{3}{4}$.
 
 In $1984$, Puglisi \cite{puglisi} made some improvements, extending further the domain of the critical strip in which it is possible to determine, up to a factor $1+o(1)$, the values of the corresponding $L$-function. \medskip\\
\textbf{Theorem. (Puglisi)} \textit{
Let $\alpha,\lambda >0$ be real numbers with $\alpha+\lambda<1$. Given $$\ell=(\log D)^{-\lambda},$$ we define the following set
$$
H(\ell, D)=\left\{s=\sigma+i t:|1-s| \geq(\log D)^{-4}, 1 / 2+\ell \leq \sigma \leq 1,|s| \leq D^{\ell / 10}\right\}
$$If $$h(-D)\leq (\log D)^{\alpha}$$then for each $s\in H(\ell,D)$ the relation $$
L(s, \chi_{D})=\frac{\zeta(2 s)}{\zeta(s)} \prod_{p \mid D}\left(1+\frac{1}{p^{s}}\right)\left[1+O\left(\exp \left\{-\frac{1}{3}(\log  D)^{1-\alpha-\lambda}\right\}\right)\right]
$$holds.
}\medskip\\An immediate consequence of Puglisi's improvement is a reformulation of Mordell's Theorem, that is, if $\zeta(\beta+i\gamma)=0$ with $\beta>\frac{1}{2}$, then, for every $\epsilon>0$, the relation $h(-D)>(\log D)^{1-\epsilon}$ holds, provided that $D>D_{0}(\beta,\gamma,\epsilon)$.\medskip

In 1934, Heilbronn \cite{hg} solved Gauss' second conjecture. He proved that, under the assumption that the general Riemann Hypothesis is not true,
$h(-D) \rightarrow \infty$ if $ D \rightarrow \infty
$. Heilbronn's result is very important, as, combined with Hecke's theorem, gives, without any assumption, that
$
h(-D) \rightarrow \infty$ if $D \rightarrow \infty
$.

In 1935, Siegel \cite{siegeltheorem} proved that $
h(-D)>D^{1 / 2-\epsilon}$ for $ D>D_{0}(\epsilon)
$ for an arbitrary $\epsilon>0$, and with a constant $D_{0}(\epsilon)$ depending only on $\epsilon$, where the constant $D_{0}(\epsilon)$ is ineffective (for alternative proofs of Siegel's Theorem see Estermann \cite{estermannsiegel}, Chowla \cite{chowla}, Goldfeld \cite{goldsiegel}, Linnik \cite{linniksiegel}, Pintz  \cite{pintzsiegel}).

Heilbronn played a fundamental role also in the attempt to prove Gauss' first conjecture. In 1934, Heilbronn and Linfoot \cite{heillin} showed that, except for the known values $-D=-3,-4,-7,-8,-11,-19,-43,-67,-163$, there is at most a tenth negative fundamental discriminant with class number one.

In 1935 Landau \cite{landau2} proved that if $
h(-D)=h
$, then the inequality $
D \leq D(h)=C h^{8} \log ^{6}(3 h)
$ holds, where $C$ is an absolute effective constant, with the possible exception of at most one negative fundamental discriminant.

In 1950 Tatuzawa \cite{tatu} proved Landau's theorem mentioned above with $
D(h)=C h^{2} \log ^{2}(13 h)
$. Furthermore, Tatuzawa made some improvements regarding the effectivization of Siegel's Theorem, showing that if $h(-D) \leq D^{1 / 2-\epsilon}$, then the inequality $D \leq D_{0}^{\prime}(\epsilon)=\max \left(e^{12}, e^{1 / \epsilon}\right)
$ holds, with the possible exception of at most one negative fundamental discriminant.\\ Finally, in 1966-1967, Baker \cite{baker1} and Stark \cite{stark1} proved independently that there is no tenth imaginary quadratic field with class number one.

The results found by Deuring \cite{dg} and Heilbronn \cite{hg} regarding the influence of the non-trivial zeros of both $\zeta(s)$ and $L(s, \chi)$ (where $\chi$ is an arbitrary real or complex character) on the real zeros of other real $L$-functions caught the interest of Linnik, who deeply analyzed this phenomenon, known as the \textit{Deuring-Heilbronn phenomenon}, in his work concerning the least prime in an arithmetic progression, finding important new results \cite{linnikheil}.\medskip\\
\textbf{Theorem. (Linnik)} 
\textit{If an $L$-function belonging to a real non-principal character modulus $D$ has a real zero $1-\delta$ with
$$
\delta \leq\frac{ A_{1}}{  \log D},
$$
then all the $L$-functions belonging to characters modulus $D$ have no zero in the domain
 $$\sigma \geq 1-\frac{A_{2}}{\log D(|t|+1)} \log \left(\frac{e A_{1}}{\delta \log D(|t|+1)}\right), \quad\quad \delta \log D(|t|+1) \leq A_{1},$$where $A_{1}$ and $A_{2}$ are absolute constants.}\medskip
 
Some improvements related to the Heilbronn phenomenon were found by Pintz in 1975 \cite{pintz4}. In particular, using elementary methods, he proved the following result.\medskip\\
 \textbf{Theorem (Pintz)} \textit{Let $L(s,\chi_{k})$ be a Dirichlet's $L$-function belonging to the non principal character (real or complex) $\chi_{k}$ modulus $k$. Suppose that  $L(s,\chi_{k})$ has a zero $s_{0}=1-\gamma +it$ with $\gamma <0.05$.\\ Then, for an arbitrary real non-principal character $\chi_{D}$ mod $D$ (for which $\chi_{k} \chi_{D}$ is also non-principal) the inequality
$$
L\left(1, \chi_{D}\right)>\frac{1}{140 U^{6 \gamma} \log ^{3} U}
$$
holds, where $U=k\left|s_{0}\right| D .$}\medskip\medskip\medskip

The aim of this work is to further investigate both the Deuring phenomenon and the Heilbronn phenomenon. We will find better estimates regarding the influence of zeros of $\zeta(s)$ on the exceptional zeros and that of the non-trivial zeros of arbitrary $L$-functions belonging to non-principal characters on the exceptional zeros, respectively.\\Regarding the Deuring phenomenon, combining elementary methods with some tools of complex analysis based on Pintz's \cite{pintz3} and Puglisi's \cite{puglisi} approach, we will go further into the critical strip. More precisely, we will prove the following theorem, provided that $L(s,\chi)$ is a Dirichlet $L$-function belonging to the real primitive character $\chi$ modulus $q$ satisfying $\chi(-1)=-1$ and $h(-q)$ is the number of classes of the imaginary quadratic field $\mathbb{Q}(\sqrt{-q})$.
\begin{teorema}\label{deuring}
Let $\eta,\ \mu>0$ be real numbers with $\eta>\max(\mu,1)$. Given $$\ell=(\log\log q)^{-\mu},$$ we define the following set
$$
H(\ell, q)=\left\{s=\sigma+i t:|1-s| \geq(\log q)^{-4}, 1 / 2+\ell \leq \sigma \leq 1,|s| \leq q^{\ell / 10}\right\}
$$If $$h(-q)\leq \frac{\log q}{(\log\log q)^{\eta}}$$then for each $s\in H(\ell,q)$ the relation $$
L(s, \chi)=\frac{\zeta(2 s)}{\zeta(s)} \prod_{p \mid q}\left(1+\frac{1}{p^{s}}\right)\left[1+O\left(\exp \left\{-\frac{1}{3}(\log \log q)^{\eta-\mu}\right\}\right)\right]
$$holds.  \end{teorema}\medskip As an immediate consequence, a new reformulation of Mordell's Theorem follows from Theorem \ref{deuring}.
\begin{co}
If $\zeta(\beta+i\gamma)=0$ with $\beta>1/2$, then for every $\eta>1$ the relation $$h(-q)>\frac{\log q}{(\log\log q)^{\eta}} $$holds, provided that $q>q_{0}(\beta,\gamma,\eta)$.
\end{co}The improvements regarding the Deuring phenomenon stated above make sense, as the inequality $h(-q)>c\log q/(\log\log q)^{\eta}$ had never been generalized to an arbitrary modulus $q$, but it was valid only for $q$ prime (\cite{gold},\cite{grosszaiger}).
\medskip\medskip\\
Regarding the Heilbronn phenomenon, we will improve Pintz's theorem stated above, showing that it is possible to extend the range of values for $\gamma$ to $0<\gamma<\frac{1}{4}$ if $\chi_{k}$ is real and to $0<\gamma\leq \frac{1}{8}$ if $\chi_{k}$ is complex. More precisely, we will use elementary methods based on Pintz's approach \cite{pintz4} to prove the following theorem.
\begin{teorema}\label{mioh1}Let $L(s,\chi_{k})$ be a Dirichlet $L$-function belonging to the real non principal character $\chi_{k}$ modulus $k$. Suppose that  $L(s,\chi_{k})$ has a zero $s_{0}=1-\gamma +it$ with $0<\gamma <\frac{1}{4}$.\\ Then, for an arbitrary real non-principal character $\chi_{D}$ mod $D$ (for which $\chi_{k} \chi_{D}$ is also non-principal) the inequality
$$L(1,\chi_{D})\geq \frac{c_{1}}{U^{b\gamma}\log^{3}U},\quad\quad \text{for}\quad \frac{1}{2(1-3\gamma)}<b<\frac{1}{2\gamma}$$
holds, where $U=k\left|s_{0}\right| D $ and $c_{1}$ is an effective constant.\\The same result can be obtained if $\chi_{k}$ is a complex non principal character, provided that $0<\gamma\leq \frac{1}{8}$.
\end{teorema}
\medskip

Theorem \ref{mioh1} has some important consequences.\\First of all, we can deduce that a zero in the half-plane $\sigma>\frac{3}{4}$ for real characters or a zero in the half-plane $\sigma\geq\frac{7}{8}$ for complex characters implies that $h(-D)\rightarrow\infty$.\\Furthermore, a weakened form of Linnik Theorem \cite{linnikheil} can be deduced (the following theorem is an improvement of Theorem 2 of \cite{pintz4}).
\begin{teorema}\label{mioh2}
If an $L$-function belonging to a non-principal character $\chi_{k}$ modulus $ k$ has a zero $s_{0}=1-\gamma+it$ with $\gamma<\frac{1}{4}$ if $\chi_{k}$ is real or $\gamma\leq\frac{1}{8}$ if $\chi_{k}$ is complex, and another $L$-function belonging to the real non-principal character $\chi_{D}$ (for which $\chi_{k} \chi_{D}$ is also non-principal) modulus $D$ has a real exceptional zero $1-\delta$, then the inequality
$$
\delta>\frac{c_{1}}{ U^{b \gamma} \log ^{5} U}\quad\quad\text{for}\quad \frac{1}{2(1-3\gamma)}<b<\frac{1}{2\gamma}
$$
holds, where $U=k\left|s_{0}\right| D$ and $c_{1}$ is the costant of Theorem \ref{mioh1}.
\end{teorema} An immediate consequence is Linnik's Theorem, stated above, in the following form.
\begin{co}
If an $L$-function belonging to a real non-principal character modulus $D$ has a real zero $1-\delta$ with
$$
\delta =O_{\epsilon}\left(\frac{1}{\log^{5+\epsilon}D}\right)\quad (\epsilon>0)
$$
then all the $L$-functions belonging to characters modulus $D$ have no zero in the domain
$$\sigma\geq 1-\frac{1}{b\log U}\log\left(\frac{c_{1}}{\delta\log^{5}U}\right)$$where $c_{1}$ and $b$ have been defined in Theorem \ref{mioh1}.
\end{co}Furthermore, from Theorem \ref{mioh2}, combined with Hecke's Theorem (see Pintz \cite{pintz1}, p. 58),  we obtain the following result regarding real zeros of real $L$-functions (the following theorem is an improvement of Theorem 3 of \cite{pintz4}).
\begin{co}
For an arbitrary $\gamma,\  0<\gamma\leq \frac{1}{8}$, there is at most one $D$, and at most one primitive real character $\chi_{D}$ modulus $D$, such that $L\left(s, \chi_{D}\right)$ vanishes somewhere in the interval
\begin{equation*}
\left[1-\min \left(\gamma, \frac{c_{1}}{  32 \log ^{5} D \cdot D^{b\gamma}}\right), 1\right]
\end{equation*}where both $c_{1}$ and $b$ have been defined in Theorem \ref{mioh1}.
\end{co}

\section{Proof of Theorem \ref{deuring}}
In order to prove Theorem \ref{deuring}, following Pintz's \cite{pintz3} and Puglisi's \cite{puglisi} approach to the Deuring-phenomenon, we need some lemmas.\\First of all, we define the function $g(n)=\sum_{d|n}\chi(d)$.
\begin{lemma}\label{perron}
Given $\frac{1}{2}+\ell\leq\sigma\leq \frac{7}{8}$ and $x\gg q $, the relation
\begin{align*}
    \sum_{n \leq x}\frac{g(n)}{n^{s}}\left(1-\frac{n}{x}\right)^{2}=&L(s, \chi) \zeta(s)+\frac{2 x^{1-s} L(1, \chi)}{(1-s)(2-s)(3-s)}+\\&+O\left(|s| \log ^{2}(2+|s|) \exp \left\{-\frac{1}{2}\frac{\log q}{(\log\log q)^{\mu}}\right\}\right)
\end{align*}
holds.
\begin{proof}Following exactly the proof of Lemma 2 of \cite{puglisi}, we obtain again that
\begin{align*}
    \frac{1}{2} \sum_{n \leq x} \frac{g(n)}{n^{s}}\left(1-\frac{n}{x}\right)^{2}&=\frac{1}{2 \pi i} \int_{-\sigma-i \infty}^{-\sigma+i \infty} \frac{L(s+w, \chi) \zeta(s+w )}{w(w+1)(w+2)} x^{w}dw+\\&\ \ \ +\frac{L(s,\chi)\zeta(s)}{2}+\frac{ x^{1-s} L(1, \chi)}{(1-s)(2-s)(3-s)}
\end{align*}
Now, using both the hypothesis of Lemma \ref{perron} and the classical estimates that were already used in the proof of Lemma 2 of \cite{puglisi}, namely$$
\zeta(i t) \ll \sqrt{|t|+1} \log (|t|+2)
$$$$
L(i t, \chi) \ll \sqrt{q(|t|+1)} \log (q(|t|+1)),
$$ we get
\begin{align*}
    &\left|\int_{-\infty}^{\infty}\frac{x^{-\sigma+iu}\sqrt{q}|t+u|\log q\log^{2}(|t+u|+2)}{(-\sigma+iu)(-\sigma+1+iu)(-\sigma+2+iu)}du\right|\ll\\&\ll |s|q^{-\ell}\log q\log^{2}(2+|s|)\ll\\&\ll|s|\log^{2}(2+|s|)\exp\left\{-\frac{\log q}{(\log\log q)^{\mu}}+\log\log q\right\}\ll\\&\ll|s|\log^{2}(2+|s|)\exp\left\{-\frac{1}{2}\frac{\log q}{(\log\log q)^{\mu}}\right\}
\end{align*}
Finally, from the above estimate we can conclude that the relation \begin{align*}
    \sum_{n \leq x}\frac{g(n)}{n^{s}}\left(1-\frac{n}{x}\right)^{2}=&L(s, \chi) \zeta(s)+\frac{2 x^{1-s} L(1, \chi)}{(1-s)(2-s)(3-s)}+\\&+O\left(|s| \log ^{2}(2+|s|) \exp \left\{-\frac{1}{2}\frac{\log q}{(\log\log q)^{\mu}}\right\}\right)
\end{align*}holds for $\frac{1}{2}+\ell\leq\sigma\leq \frac{7}{8}$ and $x\gg q $.
\end{proof}
\end{lemma}
\begin{lemma}
If $s\in H(\ell,q)$ and the following inequality $$h(-q)\leq \frac{\log q}{(\log\log q)^{\eta}}$$holds, then the relation $$
\sum_{n \leq q} \frac{g(n)}{n^{s}}=L(s, \chi) \zeta(s)+O\left(\exp \left\{-\frac{1}{3}\frac{\log q}{(\log\log q)^{\mu}}
\right\}\right)
$$holds.
\begin{proof}
First of all, we suppose that $\frac{1}{2}+\ell\leq \sigma\leq \frac{7}{8}$.\\
We know that$$\sum_{n \leq q}\frac{g(n)}{n^{s}}\left(1-\frac{n}{q}\right)^{2}=\sum_{n\leq q}\frac{g(n)}{n^{s}}-\frac{2}{q}\sum_{n\leq q}\frac{g(n)n}{n^{s}}+\frac{1}{q^{2}}\sum_{n\leq q}\frac{g(n)n^{2}}{n^{s}}$$
Using Lemma \ref{perron} we have
\begin{align*}
    \sum_{n \leq q}\frac{g(n)}{n^{s}}=&L(s, \chi) \zeta(s)+\frac{2 q^{1-s} L(1, \chi)}{(1-s)(2-s)(3-s)}+\\&+O\left(|s| \log ^{2}(1+|s|) \exp \left\{-\frac{1}{2}\frac{\log q}{(\log\log q)^{\mu}}\right\}\right)+\\&+\frac{2}{q}\sum_{n\leq q}\frac{g(n)}{n^{s-1}}-\frac{1}{q^{2}}\sum_{n\leq q}\frac{g(n)}{n^{s-2}}
\end{align*}Now, if we use Dirichlet's Class Number Formula (see Davenport \cite{dav1}, chapter $6$) and Lemma 1 of \cite{puglisi}, it follows that
\begin{align*}
    \frac{1}{q}\sum_{n\leq q}\frac{g(n)}{n^{s-1}}&\ll q^{-\frac{1}{2}-\ell}\sum_{n\leq q}g(n)\ll q^{\frac{1}{2}-\ell}L(1,\chi)\ll q^{-\ell}h(-q)\ll q^{-\ell}\frac{\log q}{(\log\log q)^{\eta}}=\\&=\exp\left\{-\ell\log q+\log\log q-\eta\log\log\log q\right\}\ll\\&\ll \exp\left\{-\frac{\log q}{(\log\log q)^{\mu}}+\log\log q\right\}\ll \exp\left\{-\frac{1}{2}\frac{\log q}{(\log\log q)^{\mu}}\right\}
\end{align*}In the same way, since $$\frac{1}{q^{2}}\sum_{n\leq q}\frac{g(n)}{n^{s-2}}\ll q^{\frac{1}{2}-\ell}L(1,\chi)$$
$$\frac{2 q^{1-s} L(1, \chi)}{(1-s)(2-s)(3-s)}\ll q^{\frac{1}{2}-\ell}L(1,\chi),$$the same estimate as before holds.\\Furthermore, we observe that
\begin{align*}
    |s| \log ^{2}(1+|s|) \exp \left\{-\frac{1}{2}\frac{\log q}{(\log\log q)^{\mu}}\right\}&\leq q^{\frac{\ell}{10}}\left(\frac{\log q}{(\log\log q)^{\mu}}\right)^{2}\exp \left\{-\frac{1}{2}\frac{\log q}{(\log\log q)^{\mu}}\right\}
\end{align*}As a consequence, combining all the previous estimates, we can conclude that
$$\sum_{n \leq q}\frac{g(n)}{n^{s}}=L(s, \chi) \zeta(s)+O\left(\exp \left\{-\frac{1}{2}\frac{\log q}{(\log\log q)^{\mu}}\right\}\left(1+q^{\frac{\ell}{10}}\left(\frac{\log q}{(\log\log q)^{\mu}}\right)^{2}\right)\right)$$On the other hand, we have
\begin{align*}
    q^{\frac{\ell}{10}}\left(\frac{\log q}{(\log\log q)^{\mu}}\right)^{2}&=\exp\left\{\frac{\ell}{10}\log q+2\log\log q-2\mu\log\log\log q\right\}\ll\exp\left\{\frac{\log q}{10(\log\log q)^{\mu}}\right\}
\end{align*}As a result, the following estimate
\begin{align*}
    &\exp \left\{-\frac{1}{2}\frac{\log q}{(\log\log q)^{\mu}}\right\}\left(q^{\frac{\ell}{10}}\left(\frac{\log q}{(\log\log q)^{\mu}}\right)^{2}\right)\ll\\&\ll\exp\left\{-\frac{1}{2}\frac{\log q}{(\log\log q)^{\mu}}+\frac{\log q}{10(\log\log q)^{\mu}}+2\log\log q\right\}\ll\exp \left\{-\frac{1}{3}\frac{\log q}{(\log\log q)^{\mu}}\right\}
\end{align*}holds.\\So, we proved the claim for $\frac{1}{2}+\ell\leq \sigma\leq \frac{7}{8}$.\\If  $\frac{7}{8}<\sigma<1$, we can conclude as in Lemma 3 of \cite{puglisi}. 
\end{proof}
\end{lemma}
Now, we define the same sets used by Pintz \cite{pintz3} and Puglisi \cite{puglisi}:$$
A_{j}=\{n \in \mathbb{N}: p \mid n \Rightarrow \chi(p)=j\} \quad(j=-1,0,1)
$$$$
R=\left\{r=b m: b \in A_{0}, m \in A_{-1}\right\}
$$
\begin{lemma}\label{lpd}
If$$
\sum_{a \in A_{1}, 1<a \leq \sqrt{q }/2} 1 \leq h(-q)
$$then$$
\chi(p)=1 \Rightarrow p>\frac{1}{2} \exp \left\{\frac{\log q}{2(h(-q)+1)}\right\}
$$
\begin{proof}
By contradiction, we suppose that$$
\chi(p)=1 \Rightarrow p\leq \frac{1}{2} \exp \left\{\frac{\log q}{2(h(-q)+1)}\right\}
$$Since $h(-q)\geq 1$, we have
$$p^{h(-q)+1}\leq \frac{1}{2^{(h(-q)+1)}}\exp\left\{\frac{\log q}{2}\right\}\leq \frac{1}{4}\sqrt{q}\leq\frac{1}{2}\sqrt{q}$$Then, we consider $p,p^{2},\dots, p^{h(-q)+1}$. Under these conditions, the sum$$
\sum_{a \in A_{1}, 1<a \leq \sqrt{q}/2} 1 
$$ has at least $h(-q)+1$ terms. Indeed, taken $a=p^{j}$ with $j=1,\dots,h(-q)+1$, we have $p|p^{j}$ and $\chi(p)=1$ by hypothesis.\\However, we have a contradiction because we got that $h(-q)+1\leq h(-q)$. 
\end{proof}
\end{lemma}
\begin{lemma}
If $\sigma\geq\frac{1}{2}+\ell$ and the inequality$$h(-q)\leq\frac{\log q}{(\log\log q)^{\eta}}$$holds, then the relation
$$\sum_{a \in A_{1}, 1<a \leq q} g(a) a^{-\sigma}\ll\exp\left\{-\frac{1}{10}(\log\log q)^{\eta}\right\}$$holds.
\begin{proof}
We know that $$
\sum_{a \in A_{1}, 1<a \leq q} g(a) a^{-\sigma} \leq \exp \left\{C \sum_{p \leq q\atop \chi(p)=1} p^{-\sigma}\right\}-1 \quad(C>0)
$$Furthermore, from Lemma \ref{lpd}, if$$
\sum_{a \in A_{1}, 1<a \leq \sqrt{q }/2} 1 \leq h(-q)
$$then$$
\chi(p)=1 \Rightarrow p>\frac{1}{2} \exp \left\{\frac{\log q}{2(h(-q)+1)}\right\}=R_{0}
$$As a result, since $\frac{1}{2}\leq \sigma<1$, $\eta>\max(\mu,1)$ and the inequalities $$1\leq h(-q)\leq \frac{\log q}{(\log\log q)^{\eta}}$$hold, we can conclude that
\begin{align*}
    &\sum_{p \leq \sqrt{q}/2\atop \chi(p)=1} p^{-\sigma}\leq 2^{\sigma}h(-q)\exp\left\{-\frac{\sigma\log q}{2(h(-q)+1)}\right\}\leq 2\frac{\log q}{(\log\log q)^{\eta}}\exp\left\{-\frac{\left(\frac{1}{2}+\ell\right)\log q}{2+2\frac{\log q}{(\log\log q)^{\eta}}}\right\}\leq\\&\leq 2\frac{\log q}{(\log\log q)^{\eta}}\exp\left\{-\frac{\left(\frac{1}{2}+\ell\right)\log q(\log\log q)^{\eta}}{4\log q}\right\}=\\&=2\frac{\log q}{(\log\log q)^{\eta}}\exp\left\{-\frac{1}{8}(\log\log q)^{\eta}-\frac{1}{2}(\log\log q)^{\eta-\mu}\right\}\leq\\&\leq2\frac{\log q}{(\log\log q)^{\eta}}\exp\left\{-\frac{1}{8}(\log\log q)^{\eta}\right\}\ll \exp\left\{-\frac{1}{10}(\log\log q)^{\eta}\right\}
\end{align*}Furthermore, for $\sigma\geq \frac{1}{2}+\ell$ we have
$$
\sum_{\sqrt{q} / 2<p \leq q,\  \chi(p)=1} p^{-\sigma} \leq \sum_{\sqrt{q} / 2<n \leq q} g(n) n^{-\sigma} \ll q^{-\frac{\ell}{2}} \sum_{\sqrt{q} / 2<n \leq q} g(n) n^{-\frac{1}{2}}
$$Even more, using Lemma A of \cite{puglisi} (for the proof see Goldfeld \cite{gold}, p. 637) with $\epsilon=\frac{1}{11}$, we have, for $0<10y<x$, $$
\sum_{y<n \leq x} \frac{g(n)}{\sqrt{n}}=\sum_{d \leq \sqrt{x}} \frac{1}{d} \sum_{y / d^{2}<k \leq x / d^{2}} \nu(k) k^{-1 / 2} \ll L(1, \chi)\left\{\frac{\sqrt{q x}}{\sqrt{y}}+\sqrt{x}+x^{\frac{1}{2}-\frac{1}{11}} q^{\frac{1}{11}}\right\}
$$Following the argument used by Puglisi in \cite{puglisi}, if we take $H=\frac{\log 4q}{\log 121}$, we obtain that
\begin{align*}
    &\sum_{\sqrt{q} / 2<p \leq q,\  \chi(p)=1} p^{-\sigma}\ll q^{-\frac{\ell}{2}} \sum_{h \leq H} \sum_{\sqrt{q} / 2\cdot(11)^{h-1}<n \leq \sqrt{q} / 2\cdot(11)^{h}} g(n) n^{-\frac{1}{2}}\ll\\&\ll q^{-\frac{\ell}{2}} L(1, \chi) \sum_{h \leq H}\left\{\sqrt{q}+\sqrt{\frac{\sqrt{q}(11)^{h}}{2}}+\left(\frac{\sqrt{q}(11)^{h}}{2}\right)^{\frac{1}{2}-\frac{1}{11}} q^{\frac{1}{11}}\right\}\ll\\&\ll H q^{\frac{1-\ell}{2}} L(1, \chi)=\frac{\log 4q}{\log 121}q^{\frac{1}{2}}q^{-\frac{1}{2(\log\log q)^{\mu}}}L(1,\chi)\ll 
    \exp\left\{-\frac{1}{8}\frac{\log q}{(\log\log q)^{\mu}}\right\}
\end{align*}where we used the estimate
$$q^{-\ell}h(-q)\ll \exp\left\{-\frac{1}{2}\frac{\log q}{(\log\log q)^{\mu}}\right\}$$ Adding both the terms, the conclusion follows.
\end{proof}
\end{lemma}
\begin{lemma}
If $\sigma\geq\frac{1}{2}+\ell$ and the inequality$$h(-q)\leq\frac{\log q}{(\log\log q)^{\eta}}$$holds, then the relation$$\sum_{n \leq q} g(n) n^{-s}=\sum_{r \in R, r\leq q} g(r) r^{- s}+O\left(\exp\left\{-\frac{1}{16}(\log\log q)^{\eta}\right\}\right)$$holds.
\begin{proof}
First of all, we observe that$$
\sum_{n \leq q} g(n) n^{-s}=\sum_{r \in R, r\leq q} g(r) r^{- s}+O\left(\sum_{r \in R, r \leq q} g(r) r^{-\sigma} \sum_{a\in A_{1},\ 1<a \leq q} g(a) a^{-\sigma}\right)
$$and$$
\sum_{r \in R, r \leq q} g(r) r^{-\sigma} \leq \sum_{b \in A_{0}} \frac{\mu^{2}(b)}{b^{\frac{1}{2}+\ell}} \sum_{k \geq 1} k^{-1-2 \ell} \ll \frac{1}{\ell} \exp \left\{\sum_{p \mid q} \frac{1}{\sqrt{p}}\right\}
$$where $\mu$ is Möbius' Function.\\Since $$h(-q)\leq\frac{\log q}{(\log\log q)^{\eta}}$$and$$
\sum_{p \mid q} 1 \leq 1+\frac{\log (h(-q))}{\log 2},
$$then$$
\exp \left\{\sum_{p \mid q} \frac{1}{\sqrt{p}}\right\}\leq\exp\left\{1+\frac{\log(h(-q))}{\log 2}\right\}\leq 3\exp\left\{\frac{\log(h(-q))}{\log 2}\right\}\ll\left(\frac{\log q}{(\log\log q)^{\eta}}\right)^{\frac{1}{\log 2}}
$$It follows that$$
\sum_{r \in R, r \leq q} g(r) r^{-\sigma}\ll\frac{1}{\ell}\left(\frac{\log q}{(\log\log q)^{\eta}}\right)^{\frac{1}{\log 2}}=(\log\log q)^{\mu-\frac{\eta}{\log 2}}(\log q)^{\frac{1}{\log 2}}$$As a result, we have
\begin{align*}
 &\sum_{r \in R, r \leq q} g(r) r^{-\sigma} \sum_{a\in A_{1},\ 1<a \leq q} g(a) a^{-\sigma}\ll\\&\ll (\log\log q)^{\mu-\frac{\eta}{2}}(\log q)^{\frac{1}{\log 2}}\left(\exp\left\{-\frac{1}{10}(\log\log q)^{\eta}\right\}\right)\ll \exp\left\{-\frac{1}{16}(\log\log q)^{\eta}\right\}
\end{align*}
\end{proof}
\end{lemma}
\begin{lemma}If $\sigma\geq\frac{1}{2}+\ell$ and the inequality$$h(-q)\leq\frac{\log q}{(\log\log q)^{\eta}}$$holds, then the relation
\begin{align*}
    \sum_{r \in R \atop r \leq q} \frac{g(r)}{r^{s}}=&\zeta(2 s) \prod_{p\mid q}\left(1+\frac{1}{p^{s}}\right)\left[1+O\left(\exp \left\{-\frac{1}{3}(\log \log q)^{\eta-\mu}\right\}\right)\right]+\\&+O\left(\exp\left\{-\frac{1}{2}\frac{\log q}{(\log\log q)^{\mu}}\right\}\right)
\end{align*}
holds.
\begin{proof}
We have already seen that$$
\frac{1}{\ell} \sum_{h \mid q} \frac{\mu^{2}(h)}{\sqrt{h}}=\frac{1}{\ell} \prod_{p \mid q}\left(1+\frac{1}{\sqrt{p}}\right) \ll (\log\log q)^{\mu-\frac{\eta}{\log 2}}(\log q)^{\frac{1}{\log 2}}
$$Furthermore, if $n\not\in R$, then $n>R_{0}$.\\It follows that
\begin{align*}
 &\sum_{r \in R \atop r \leq q} \frac{g(r)}{r^{s}}=\sum_{h|q}\frac{\mu^{2}(h)}{h^{s}}\sum_{r\in R,\ r\leq \sqrt{q/h}}r^{-2s}=\\&=\sum_{h|q}\frac{\mu^{2}(h)}{h^{s}}\left[\zeta(2s)+O\left(\sum_{r>R_{0}}r^{-1-2\ell}\right)\right]+O\left(\sum_{h|q}\frac{\mu^{2}(h)}{h^{\frac{1}{2}+\ell}}\sum_{r>\sqrt{q/h}}r^{-1-2\ell}\right)=\\&=\prod_{p \mid q}\left(1+\frac{1}{p^{s}}\right)\left[\zeta(2s)+O\left(\frac{1}{\ell}\exp\left\{-\frac{\log q}{2(\log\log q)^{\mu}h(-q)}\right\}\right)\right]+\\&\ \ \  +O\left(\frac{q^{-\ell}}{\ell}\sum_{h|q}\frac{\mu^{2}(h)}{\sqrt{h}}\right)=\\&=\prod_{p \mid q}\left(1+\frac{1}{p^{s}}\right)\left[\zeta(2s)+O\left(\frac{1}{\ell}\exp\left\{-\frac{1}{2}(\log\log q)^{\eta-\mu}\right\}\right)\right]+\\&\ \ \  +O\left(q^{-\ell}(\log\log q)^{\mu-\frac{\eta}{\log 2}}(\log q)^{\frac{1}{\log 2}}\right)=\\&=\prod_{p \mid q}\left(1+\frac{1}{p^{s}}\right)\left[\zeta(2s)+O\left((\log\log q)^{\mu}\exp\left\{-\frac{1}{2}(\log\log q)^{\eta-\mu}\right\}\right)\right]+\\&\ \ \  +O\left(\exp\left\{-\frac{\log q}{(\log\log q)^{\mu}}+\frac{1}{\log 2}\log\log q+\left(\mu-\frac{\eta}{\log 2}\right)\log\log\log q\right\}\right)=\\&=\zeta(2 s) \prod_{p\mid q}\left(1+\frac{1}{p^{s}}\right)\left[1+O\left(\exp \left\{-\frac{1}{3}(\log \log q)^{\eta-\mu}\right\}\right)\right]+\\&\ \ \  +O\left(\exp\left\{-\frac{1}{2}\frac{\log q}{(\log\log q)^{\mu}}\right\}\right)
\end{align*}
\end{proof}
\end{lemma}
Now, we are ready to prove Theorem \ref{deuring}.\\
Using all the results we found previously, we can conclude that
\begin{align*}
    &L(s,\chi)\zeta(s)=\sum_{n \leq q} \frac{g(n)}{n^{s}}+O\left(\exp \left\{-\frac{1}{3}\frac{\log q}{(\log\log q)^{\mu}}
\right\}\right)=\\&=\sum_{r \in R, r\leq q} g(r) r^{- s}+O\left(\exp\left\{-\frac{1}{16}(\log\log q)^{\eta}\right\}\right)=\\&=\zeta(2 s) \prod_{p\mid q}\left(1+\frac{1}{p^{s}}\right)\left[1+O\left(\exp \left\{-\frac{1}{3}(\log \log q)^{\eta-\mu}\right\}\right)\right]+\\&\ \ \ \ +O\left(\exp\left\{-\frac{1}{16}(\log\log q)^{\eta}\right\}\right)=\\&=\zeta(2 s) \prod_{p\mid q}\left(1+\frac{1}{p^{s}}\right)\left[1+O\left(\exp \left\{-\frac{1}{3}(\log \log q)^{\eta-\mu}\right\}\right)\right]
\end{align*}
\section{Proof of Theorem \ref{mioh1}}
Following exactly Pintz's proof of Theorem 1 of \cite{pintz4}, we define the following sets$$
A_{\nu}=\left\{n\in\mathbb{N} ;\  p|n, p \operatorname{prime} \rightarrow \chi_{D}(p)=\nu\right\} \quad(\nu=-1,0,1)
$$$$
C=\left\{c ; c=uv, u \in A_{1}, v \in A_{0}\right\}
$$and the following two multiplicative functions $$
g_{\lambda}(n)=\sum_{d \mid n} \lambda(d)=\left\{\begin{array}{lll}
1, & \text { if } & n=l^{2} \\
0, & \text { if } & n \neq l^{2}
\end{array}\right.
$$(where $\lambda(n)$ denotes Liouville's $\lambda$-function) and\begin{equation}\label{mh}
    g_{D}(n)=\sum_{d \mid n} \chi_{D}(d)=\prod_{p^{\alpha} \mid\mid n}\left(1+\chi_{D}( p)+\ldots+\chi_{D}^{\alpha}(p)\right) \geq 0
\end{equation}
Again, from Pintz's proof of Theorem 1 in \cite{pintz4}, for $n=uv m=c m,\quad$ $ u \in A_{1}, v \in A_{0}, m \in A_{-1},$ we get\begin{equation}\label{mh1}
  g_{\lambda}(n)=g_{\lambda}(u) g_{\lambda}(v) g_{\lambda}(m)=\sum_{c_{l} \mid c,\  c_{l}=u_{l} v_{l}\atop u_{l}\in A_{1},\ v_{l}\in A_{0}} 2^{\omega\left(u_{l}\right)} \lambda\left(c_{l}\right) g_{D}\left(\frac{n}{c_{l}}\right) 
\end{equation}
(where $\omega(n)$ denotes the number of distinct prime divisors of $n$) and, for $c \in C,\  c=uv,\  u \in A_{1}, \ v \in A_{0}$, we have 
\begin{equation}\label{mh2}
    2^{\omega(u)} \leq g_{D}(c) \leq d(c)
\end{equation}
Now, let $b, h$ two positive real numbers, with $1<h<2b$. Thus, considering (\ref{mh}), (\ref{mh1}) and (\ref{mh2}) we have
\begin{equation}\label{emh1}
\begin{aligned}
     &\left|\sum_{n\leq U^{b}\atop n=l^{2}}\frac{\chi_{k}(n)}{n^{s_{0}}}\right|=\left|\sum_{n \leq U^{b}} \frac{\chi_{k}(n)}{n^{s_{0}}} g_{\lambda}(n)\right|=\left|\sum_{n \leq U^{b}} \frac{\chi_{k}(n)}{n^{s_{0}}}\sum_{c\in C,\ c|n\atop c=uv,\ u\in A_{1},v\in A_{0}}2^{\omega(u)} \lambda(c) g_{D}\left(\frac{n}{c}\right)\right|=\\&=\left|\sum_{\substack{c \leqslant U^{b}, c \in C \\ c=uv, u \in A_{1}, v \in A_{0}}} \frac{2^{\omega(u)} \lambda(c) \chi_{k}(c)}{c^{s_{0}}} \sum_{r \leqslant U^{b} / c} \frac{\chi_{k}(r)}{r^{s_{0}}} g_{D}(r)\right|\leq\\&\leq \sum_{n \leq U^{b/h}} \frac{d(n)}{n^{1-\gamma}}\left|\sum_{r \leq U^{b} / n} \frac{\chi_{k}(r)}{r^{s_{0}}} g_{D}(r)\right|+\sum_{U^{b/h}<n \leq U^{b}} \frac{g_{D}(n)}{n^{1-\gamma}} \sum_{r \leq U^{b} / n} \frac{d(r)}{r^{1-\gamma}}=\\&=\sum_{1}+\sum_{2}\end{aligned}
    \end{equation}
    Before trying to estimate both the two sums in (\ref{emh1}), we find a lower bound for $$\left|\sum_{n\leq U^{b}\atop n=l^{2}}\frac{\chi_{k}(n)}{n^{s_{0}}}\right|$$
       In order to do this, we consider the two cases ($\chi_{k}$ real or complex) separately.\\We start with $\chi_{k}$ real and non principal. We observe that  \begin{equation*}
        \begin{aligned}
   \sum_{n=1\atop n=l^{2}}^{\infty}\frac{\chi_{k}(n)}{n^{s_{0}}}&=\sum_{l=1\atop (l,k)=1}^{\infty}\frac{1}{l^{2s_{0}}}=\\&=\sum_{l=1}^{\infty}\frac{\chi_{0,k}(l)}{l^{2s_{0}}}=\\&=L(2s_{0},\chi_{0,k})=\\&=\prod_{p\nmid k}\left(1-\frac{1}{p^{2s_{0}}}\right)^{-1}
        \end{aligned}
    \end{equation*}where in the last inequality we used Euler's identity.\\Hence,     \begin{equation}\label{reals}
        \begin{aligned}
\left|\sum_{n=1\atop n=l^{2}}^{\infty}\frac{\chi_{k}(n)}{n^{s_{0}}}\right|&=\left|\prod_{p\nmid k}\left(1-\frac{1}{p^{2s_{0}}}\right)^{-1}\right|\geq\\&\geq \prod_{p\nmid k}\frac{1}{1+\frac{1}{p^{2(1-\gamma)}}}\geq\\&\geq \prod_{p}\frac{1}{1+\frac{1}{p^{2(1-\gamma)}}}=\\&=\frac{\zeta(4(1-\gamma))}{\zeta(2(1-\gamma))}>\\&>\frac{\zeta(4)}{\zeta(\frac{3}{2})}>\frac{\pi^{4}}{270}>0.36
        \end{aligned}
    \end{equation} since $0<\gamma<\frac{1}{4}$.\\Now, we turn to the case where $\chi_{k}$ is complex and non principal. Since $0<\gamma\leq \frac{1}{8}$, we observe that
    \begin{equation}\label{complex}
        \begin{aligned}
       \left|\sum_{n=1\atop n=l^{2}}^{\infty}\frac{\chi_{k}(n)}{n^{s_{0}}}\right|&\geq  1-\sum_{l=2}^{\infty} \frac{1}{l^{2(1-\gamma)}}\geq\\&\geq1-\sum_{l=2}^{10}\frac{1}{l^{7/4}}-\int_{10}^{\infty}\frac{dl}{l^{7/4}}\geq 0.029
        \end{aligned}
    \end{equation}
    Now, we separately estimate the two sums of (\ref{emh1}).\\We begin with the first one: $$\sum_{1}=\sum_{n \leq U^{b/h}} \frac{d(n)}{n^{1-\gamma}}\left|\sum_{r \leq U^{b} / n} \frac{\chi_{k}(r)}{r^{s_{0}}} g_{D}(r)\right|$$We start by considering the inner sum, that is$$\left|\sum_{r \leq U^{b} / n} \frac{\chi_{k}(r)}{r^{s_{0}}} g_{D}(r)\right|$$Let $y\geq U^{b}/n$ be a fixed number and let $z$ be a parameter we will choose later.\\ Since $U=kD|s_{0}|$, we have\medskip
    \begin{equation}\label{miosum1}
        \begin{aligned}
        &\left|\sum_{r \leq y} \frac{\chi_{k}(r)}{r^{s_{0}}} g_{D}(r)\right|\leq\\&\leq \left|\sum_{d\leq z}\frac{\chi_{k}(d)\chi_{D}(d)}{d^{s_{0}}}\cdot \sum_{l\leq y/d}\frac{\chi_{k}(l)}{l^{s_{0}}}\right|+\left|\sum_{l\leq y/z}\frac{\chi_{k}(l)}{l^{s_{0}}}\cdot \sum_{z<d\leq y/l}\frac{\chi_{k}(d)\chi_{D}(d)}{d^{s_{0}}}\right|\leq\\&\leq \sum_{d\leq z}\frac{1}{d^{1-\gamma}}\cdot \frac{2|s_{0}|\sqrt{k}\log k}{\left(\frac{y}{d}\right)^{1-\gamma}}+\sum_{l\leq y/z}\frac{1}{l^{1-\gamma}}\cdot \frac{2|s_{0}|\sqrt{kD}\log (kD)}{z^{1-\gamma}}\leq\\&\leq \frac{ z\cdot 2|s_{0}|\sqrt{k}\log k}{y^{1-\gamma}}+2|s_{0}|\sqrt{kD}\log(kD)\cdot\sum_{l\leq y/z}\frac{1}{l^{1-\gamma}z^{1-\gamma}}\leq\\&\leq \frac{ z\cdot 2|s_{0}|\sqrt{k}\log k}{y^{1-\gamma}}+2|s_{0}|\sqrt{kD}\log(kD)\cdot\frac{y^{\gamma}\log\left(\frac{y}{z}\right)}{z}
        \end{aligned}
    \end{equation}where in the second step we used the Polya-Vinogradov inequality, while in the last inequality we used the fact that
    \begin{align*}
        \sum_{l\leq y/z}\frac{1}{l^{1-\gamma}z^{1-\gamma}}=\frac{y^{\gamma}}{z} \sum_{l \leq y / z} \frac{1}{l^{1-\gamma}\left(\frac{y}{z}\right)^{\gamma}}<\frac{y^{\gamma}}{z} \sum_{l \leq y / z} \frac{1}{l^{1-\gamma}l^{\gamma}}\leq\frac{y^{\gamma}\log\left(\frac{y}{z}\right)}{z}
    \end{align*}
    Now, we choose $z$ such that
    $$\frac{z}{y^{1-\gamma}}=\frac{\sqrt{D}y^{\gamma}}{z}$$or equivalently,$$z=y^{\frac{1}{2}}D^{\frac{1}{4}}$$Using this value for $z$, the relation (\ref{miosum1}) becomes
    \begin{equation*}
       \left|\sum_{r \leq y} \frac{\chi_{k}(r)}{r^{s_{0}}} g_{D}(r)\right| \leq 2\cdot y^{\gamma-\frac{1}{2}}D^{\frac{1}{4}}|s_{0}|\sqrt{k}\log k+2|s_{0}|\sqrt{k}\cdot D^{\frac{1}{4}}\cdot y^{\gamma-\frac{1}{2}}\cdot\log(kD)\log\left(\frac{\sqrt{y}}{D^{1/4}}\right)
    \end{equation*}
    Now, we consider $y=U^{b}/n$. It follows that
    \begin{equation*}
        \begin{aligned}
        \sum_{1}&=\sum_{n \leq U^{b/h}} \frac{d(n)}{n^{1-\gamma}}\left|\sum_{r \leq U^{b} / n} \frac{\chi_{k}(r)}{r^{s_{0}}} g_{D}(r)\right|\ll\\&\ll \sum_{n \leq U^{b/h}} \frac{d(n)}{n^{1-\gamma}}\cdot \left(\frac{U^{b}}{n}\right)^{\gamma-\frac{1}{2}} D^{\frac{1}{4}}\left(2|s_{0}|\sqrt{k}\log k+2|s_{0}|\sqrt{k}\log(kD)\log\left(\frac{\sqrt{y}}{D^{1/4}}\right)\right)\ll\\&\ll 2|s_{0}|\sqrt{k}\cdot  U^{b\left(\gamma-\frac{1}{2}\right)+\frac{1}{4}}\cdot \log^{2}U\cdot\sum_{n\leq U^{b/h}}\frac{d(n)}{\sqrt{n}}\ll\\&\ll 2|s_{0}|\sqrt{k}\cdot  U^{b\left(\gamma-\frac{1}{2}\right)+\frac{1}{4}}\cdot \log^{3}U\cdot \left(U^{\frac{b}{h}}\right)^{\frac{1}{2}}=\\&=2|s_{0}|\sqrt{k}\cdot  U^{b\left(\gamma-\frac{1}{2}+\frac{1}{2h}\right)+\frac{1}{4}}\cdot \log^{3}U
        \end{aligned}
    \end{equation*}
    At this point, we observe that$$b\left(\gamma-\frac{1}{2}+\frac{1}{2h}\right)+\frac{1}{4}<0$$if and only if
    \begin{equation}\label{condition}
b>\frac{1}{4}\cdot \frac{1}{\left(\frac{1}{2}-\gamma-\frac{1}{2h}\right)}\quad\quad\text{and}\quad\quad
\frac{1}{2}-\gamma-\frac{1}{2h}>0
\end{equation}Under these conditions we can conclude that the estimate
\begin{equation}\label{sumh1}
    \sum_{1}\ll 2|s_{0}|\sqrt{k}\cdot  U^{b\left(\gamma-\frac{1}{2}+\frac{1}{2h}\right)+\frac{1}{4}}\cdot \log^{3}U
\end{equation}is non trivial,
if $U\geq U_{0}(\gamma)$, where $U_{0}(\gamma)$ is a constant depending on $\gamma$.\medskip

    Now, we turn our attention to the second sum of (\ref{emh1}), that is
    $$\sum_{2}=\sum_{U^{b/h}<n \leq U^{b}} \frac{g_{D}(n)}{n^{1-\gamma}} \sum_{r \leq U^{b} / n} \frac{d(r)}{r^{1-\gamma}}$$Since
    $$\sum_{r\leq U^{b}/n}\frac{d(r)}{r^{1-\gamma}}\ll\left(\frac{U^{b}}{n}\right)^{\gamma}\sum_{r\leq U^{b}/n}\frac{d(r)}{r}\ll \left(\frac{U^{b}}{n}\right)^{\gamma}\cdot\left(\frac{1}{2}+o(1)\right)\log^{2} U$$we have
    $$\sum_{2}\ll U^{b\gamma}\cdot\log^{2}U\cdot \left(\frac{1}{2}+o(1)\right) \sum_{U^{b/h}<n \leq U^{b}} \frac{g_{D}(n)}{n}$$
However, from Lemma 1 of \cite{pintz2}, we know that
    \begin{equation*}
        \begin{aligned}
        \sum_{U^{b/h}<n \leq U^{b}} \frac{g_{D}(n)}{n}&=b\left(1-\frac{1}{h}\right)\log U\cdot L(1,\chi_{D})+O\left(\sqrt{\frac{\sqrt{D}\log D\log U}{U^{\frac{b}{h}}}}\right)=\\&=b\left(1-\frac{1}{h}\right)\log U\cdot L(1,\chi_{D})+O\left(U^{-\left(\frac{b}{2h}-\frac{1}{4}\right)}\log U\right)=\\&=\log U\cdot\left( b\left(1-\frac{1}{h}\right)L(1,\chi_{D})+O\left(U^{-\left(\frac{b}{2h}-\frac{1}{4}\right)}\right)\right)
        \end{aligned}
    \end{equation*}which is well defined, as we supposed that $1<h<2b$.\\Hence, we can conclude that
    \begin{equation}\label{sumh2}
    \begin{aligned}
    \sum_{2}&\ll U^{b\gamma}\cdot\log^{2}U\cdot \left(\frac{1}{2}+o(1)\right) \sum_{U^{b/h}<n \leq U^{b}} \frac{g_{D}(n)}{n}\ll\\&\ll U^{b\gamma}\cdot\log^{2}U\cdot \left(\frac{1}{2}+o(1)\right)\cdot \log U\cdot\left( b\left(1-\frac{1}{h}\right)L(1,\chi_{D})+O\left(U^{-\left(\frac{b}{2h}-\frac{1}{4}\right)}\right)\right)\leq \\&\leq c_{0} U^{b\gamma}\log^{3}U\cdot L(1,\chi_{D})
    \end{aligned}
    \end{equation}if $U\geq U^{\prime}_{0}(\gamma)$, where $U^{\prime}_{0}(\gamma)$ is a constant depending on $\gamma$ and $c_{0}$ is an effective constant.\medskip
    
    At this point, if $\chi_{k}$ is a real character, combining (\ref{reals}), (\ref{sumh1}), (\ref{sumh2}), under the conditions (\ref{condition}) seen above, we get
    \begin{equation*}
       0.36\leq 2|s_{0}|\sqrt{k}\cdot  U^{b\left(\gamma-\frac{1}{2}+\frac{1}{2h}\right)+\frac{1}{4}}\cdot \log^{3}U+c_{0} U^{b\gamma}\log^{3}U\cdot L(1,\chi_{D})
    \end{equation*}or equivalently,
    \begin{equation*}
        0.36- 2|s_{0}|\sqrt{k}\cdot  U^{b\left(\gamma-\frac{1}{2}+\frac{1}{2h}\right)+\frac{1}{4}}\cdot \log^{3}U\leq c_{0} U^{b\gamma}\log^{3}U\cdot L(1,\chi_{D})
    \end{equation*}
    Furthermore, for $U\geq U_{0}(\gamma)$ sufficiently large, and so $D\geq D_{0}(\gamma)$ sufficiently large, we have$$2|s_{0}|\sqrt{k}\cdot  U^{b\left(\gamma-\frac{1}{2}+\frac{1}{2h}\right)+\frac{1}{4}}\cdot \log^{3}U\leq \frac{3}{10}$$Hence, if $\chi_{k}$ is a real non principal character, we can conclude that
    $$L(1,\chi_{D})\geq \frac{1}{c_{0}U^{b\gamma}\log^{3}U}\geq \frac{c_{1}}{U^{b\gamma}\log^{3}U}$$where $c_{1}$ is an effective constant.\\ On the other hand, if $\chi_{k}$ is a complex character, combining (\ref{complex}), (\ref{sumh1}), (\ref{sumh2}), under the conditions (\ref{condition}) seen above, we get
    \begin{equation*}
       0.029\leq 2|s_{0}|\sqrt{k}\cdot  U^{b\left(\gamma-\frac{1}{2}+\frac{1}{2h}\right)+\frac{1}{4}}\cdot \log^{3}U+c_{0} U^{b\gamma}\log^{3}U\cdot L(1,\chi_{D})
    \end{equation*}or equivalently,
    \begin{equation*}
        0.029- 2|s_{0}|\sqrt{k}\cdot  U^{b\left(\gamma-\frac{1}{2}+\frac{1}{2h}\right)+\frac{1}{4}}\cdot \log^{3}U\leq c_{0} U^{b\gamma}\log^{3}U\cdot L(1,\chi_{D})
    \end{equation*}
    Furthermore, for $U\geq U_{0}(\gamma)$ sufficiently large, and so $D\geq D_{0}(\gamma)$ sufficiently large, we have$$2|s_{0}|\sqrt{k}\cdot  U^{b\left(\gamma-\frac{1}{2}+\frac{1}{2h}\right)+\frac{1}{4}}\cdot \log^{3}U\leq \frac{1}{50}$$Hence, if $\chi_{k}$ is a complex non principal character, we can conclude that
    $$L(1,\chi_{D})\geq \frac{1}{c^{\prime}_{0}U^{b\gamma}\log^{3}U}\geq \frac{c^{\prime}_{1}}{U^{b\gamma}\log^{3}U}$$where $c^{\prime}_{1}$ is an effective constant.\\Finally, we observe that, in order to have a non trivial estimate, $b$ shall satisfy $b<\frac{1}{2\gamma}$. However, due to conditions $(\ref{condition})$, we already know that$$b>\frac{1}{4}\cdot \frac{1}{\left(\frac{1}{2}-\gamma-\frac{1}{2h}\right)}$$and $$\frac{1}{2}-\gamma-\frac{1}{2h}>0$$or equivalently, $$\gamma<\frac{1}{2}-\frac{1}{2h}$$where $1<h<2b$.\\Hence, we shall have
    $$\frac{1}{4}\cdot \frac{1}{\left(\frac{1}{2}-\gamma-\frac{1}{2h}\right)}<\frac{1}{2\gamma}$$or equivalently, $$\gamma<\frac{1}{3}-\frac{1}{3h}$$Now, we observe that, for $h>1$, the inequality $$\frac{1}{3}-\frac{1}{3h}<\frac{1}{2}-\frac{1}{2h}$$ is always satisfied. As a result, provided that $h>1$ as we supposed before, $b$, $\gamma$ and $h$ shall satisfy simultaneously only the following three relations:
    \begin{equation}\label{eqmio1}
       b<\frac{1}{2\gamma}
    \end{equation}
    \begin{equation}\label{eqmio2}
        1<h<2b
    \end{equation}
    \begin{equation}\label{eqmio3}
       \gamma<\frac{1}{3}-\frac{1}{3h}
    \end{equation}
    Now, we observe that, from (\ref{eqmio1}) and (\ref{eqmio2}), the inequality $$h<\frac{1}{\gamma}$$holds.\\On the other hand, from (\ref{eqmio3}) we have$$h>\frac{1}{1-3\gamma}$$
    As a result, we get \begin{equation}\label{eqmio4}
    \frac{1}{\gamma}>h>\frac{1}{1-3\gamma}\end{equation}or even better,
    \begin{equation}\label{eqmio5}
        \frac{1}{\gamma}>2b>h>\frac{1}{1-3\gamma}
    \end{equation} From (\ref{eqmio4}) it follows that \begin{equation*}\gamma<\frac{1}{4},\end{equation*}which makes sense, since it is stated in the hypotheses for the real case, while $\gamma\leq \frac{1}{8}<\frac{1}{4}$ for the complex case.\\On the other hand, (\ref{eqmio5}) implies that \begin{equation*}\frac{1}{2(1-3\gamma)}<b<\frac{1}{2\gamma}\end{equation*} Hence, having fixed $b$ such that
    \begin{equation*}\frac{1}{2(1-3\gamma)}<b<\frac{1}{2\gamma},\end{equation*}if we choose $h$ such that
    $$\frac{1}{1-3\gamma}<h<2b,$$we have the inequality
$$L(1,\chi_{D})\geq \frac{c_{1}}{U^{b\gamma}\log^{3}U}$$if $\chi_{k}$ is real, or$$L(1,\chi_{D})\geq \frac{c^{\prime}_{1}}{U^{b\gamma}\log^{3}U}$$if $\chi_{k}$ is complex,
 where $U=k\left|s_{0}\right| D $ and $c_{1}, c_{1}^{\prime}$ are effective constants.\\The proof of Theorem \ref{mioh1} is complete.
 \section{Proof of Theorem \ref{mioh2}}
 As in the proof of Theorem 2 of \cite{pintz4}, by a result of Page \cite{page}, given $\chi_{D}$ a real non-principal character $\bmod D,$ we know that the greatest real zero $1-\delta$ of $L\left(s, \chi_{D}\right)$ satisfies
$$
\frac{L(1,\chi_{D})}{\delta} \leq \log ^{2} D 
$$Furthermore, since $U=k|s_{0}|D$ by hypothesis, then $\log^{2}D\leq \log ^{2}U.$  Hence,$$
\frac{L(1,\chi_{D})}{\delta} \leq \log ^{2} U
$$Now, using Theorem \ref{mioh1}, it follows that$$
\delta>\frac{c_{1}}{ U^{b \gamma} \log ^{5} U}\quad\text{for }\frac{1}{2(1-3\gamma)}<b<\frac{1}{2\gamma}.$$
\medskip\medskip\medskip\medskip\\
\addcontentsline{toc}{chapter}{\bibname}

\end{document}